\documentclass{ifacconf}

\usepackage{graphicx}      
\usepackage{subfig}
\usepackage{natbib}
\usepackage{amsmath,amssymb,amsopn,cases,mathrsfs,stmaryrd}
\usepackage{hhline}
\usepackage{blindtext,xcolor}
\usepackage{multirow}
\usepackage[bb=boondox]{mathalfa}

\newcommand{\R}{\mathbb{R}}                         
\renewcommand{\H}{\mathcal{H}}                      
\renewcommand{\L}{\mathcal{L}}                      
\renewcommand{\O}{\Omega}                           
\newcommand{\J}{\mathcal{J}}                        
\newcommand{\pardev}[2]                             
{\ensuremath{\frac{\partial #1}{\partial #2}}}

\DeclareMathOperator*{\argmin}{argmin}              

\begin{document}

\begin{frontmatter}

\title{Learning-based Traffic State Reconstruction using Probe Vehicles\thanksref{footnoteinfo}} 

\thanks[footnoteinfo]{This research is partially funded by the KAUST Office of Sponsored Research under Award No. OSR-2019-CRG8-4033, the Swedish Foundation for Strategic Research and Knut and Alice Wallenberg Foundation. The authors are affiliated with the Wallenberg AI, Autonomous Systems and Software Program (WASP).}

\author[KTH]{John Liu} 
\author[KTH,corresponding]{Matthieu Barreau} 
\author[KTH]{Mladen \v{C}i\v{c}i{\'c}}
\author[KTH]{Karl H. Johansson}

\address[KTH]{Division of Decision and Control Systems, KTH Royal Institute of Technology Stockholm, Sweden (e-mail: barreau,johnliu,cicic,kallej@kth.se).}
\address[corresponding]{Corresponding author}

\begin{abstract}                
This article investigates the use of a model-based neural network for the traffic reconstruction problem using noisy measurements coming from Probe Vehicles (PV). The traffic state is assumed to be the density only, modeled by a partial differential equation. There exist various methods for reconstructing the density in that case. However, none of them perform well with noise and very few deal with lagrangian measurements. This paper introduces a method that can reduce the processes of identification, reconstruction, prediction, and noise rejection into a single optimization problem. Numerical simulations, based either on a macroscopic or a microscopic model, show good performance for a moderate computational burden.
\end{abstract}


\end{frontmatter}

\maketitle

\section{Introduction}

Reducing traffic emissions while keeping the travel time low is undoubtedly a major concern \citep{ferrara2018}. Controlling the traffic is one solution that does not imply a modification of the road infrastructure. Nevertheless, the current control laws \citep{PIACENTINI201813,cicic2019multiclass} require to estimate the traffic density. 

To reconstruct the density, one can derive a backstepping observer \citep{YU2019183}.
However, this does not allow the observation of mixed regimes (both free and congested roads) which is the critical phase for traffic control. Dynamic boundary observers as the one by \citet{qi2018boundary} can correct this issue. However, the state reconstruction is poor and it requires a huge amount of measurements.

One solution to get more measurements without affecting the traffic infrastructure is to use data from PV \citep{herrera2010incorporation}. This is possible thanks to GPS measurements \citep{amin2008mobile} or probing vehicles \citep{seo2015probe}. An empirical reconstruction using PV has been done by \cite{cicic2020numerical}. It is simple and cannot effectively deal with sparse or noisy measurements. Unfortunately, this is often the case in real-world situations when the capacity for data collection may be restrictive.

Instead, one can perform density estimation using a traffic model. The estimation problem then boils down to solving this system of equations. This was solved mathematically by \cite{delle2019} where the authors developed a numerical scheme to estimate the bulk traffic density. These results were further extended by \cite{barreau2020} where the problem was considered from a control perspective. The obtained observer uses a similar model for traffic flow with an exponential or finite-time convergence of the error. 
These methods indeed solve the traffic flow problem with noiseless data. However, they encounter complications when dealing with inaccurate measurements, even for the simplest traffic flow model.

\cite{raissi2019} present an alternative model-based approach using a physics informed neural network to solve Partial Differential Equations (PDE). 
The authors introduce a mechanism into the neural network such that it is capable of solving a PDE. From a control perspective, this can be seen as a model-based neural network and its usefulness for traffic state reconstruction has been discussed by \citep{huang2020} for instance. We intend to further investigate this problem with a special focus on noise rejection, prediction, and concrete application using a high-fidelity microscopic traffic simulator.

In this paper, we use a physics informed neural network to estimate the density in two traffic scenari using measurements from PV. In the first case, we have full knowledge of the PDE which defines the physics of the simulation while in the second one, we aim at identifying the appropriate model for the system to estimate the traffic density. Lastly, the contribution of this paper lies in that we use a neural network to reduce the problem of identification, reconstruction, prediction, and noise rejection into one single optimization problem.

In Section \ref{section:traffic_flow}, we discuss several micro- and macroscopic models for traffic flow and then present a coupled micro-macro model as well as formulate the general problem that we intend to solve. In Section \ref{section:machine_learning}, we describe the fundamentals of learning. We then present the reconstruction results for a finite difference scheme and an external microscopic simulation in Section \ref{section:results} and, lastly, offer a few concluding remarks in Section \ref{section:conclusion}.

\section{Traffic Flow Theory} \label{section:traffic_flow}

There are two main approaches to model traffic flow, namely microscopic and macroscopic \citep{ferrara2018}. A microscopic model defines the dynamic of each vehicle in a traffic stream. This enables us to observe specific phenomena that are exclusively due to the behavior of individual drivers. The main drawbacks are the high computational demand and large parameter space.

On the other hand, a macroscopic model defines the dynamics employing aggregate variables in a traffic flow system. This reduces the computational load and model complexity. The compromise is that a macroscopic model can only simulate general behaviors in traffic flow rather than individual interactions between drivers.

In this section, we present first-order models for microscopic and macroscopic simulations. We then combine them into a coupled micro-macro model which we use to define the problem formulation.

\subsection{Microscopic Model}

A common model for microscopic simulation is the first order Follow-the-Feader (FL1) model which was rigorously formulated by \cite{argall2002}:
\begin{equation} \label{equ:micro1}
    \begin{cases}
        \dot{x}_i(t) = V(\rho_i(t)), & i=1,2,\dots,n, \\
        \dot{x}_{n+1}(t) = V_L(t), \\
    \end{cases}
\end{equation}
with appropriate initial positions. It consists of one leading vehicle at position $x_{n+1}$ driving with speed $V_L(t)$ and $n$ following vehicles at positions $x_1$, $x_2$, \dots, $x_n$. The dynamics of the following vehicles are described by the velocity function $V \in C ([0,1], [0,V_{max}])$, where
\begin{equation} \label{equ:micro_density}
    \rho_i(t) = \frac{l_n}{x_{i+1}(t)-x_i(t)} \in [0,1]
\end{equation}
is the normalized density between $x_{i+1}$ and $x_i$. To avoid collisions $x_{i+1} - x_i > l_n$ must hold. It is satisfied through the inverse dependency between the relative density and the difference in position as long as the inequality is true for the initial positions. Moreover, $l_n=\frac{m}{n}$ is the average length of the vehicles and $m$ is the total mass of the all vehicles. The velocity of the leading vehicle $V_L(t)$ can be chosen freely such that $V_L\geq 0$, to avoid non-physical behaviors. Furthermore, the velocity function $V(\rho)$ needs to be a bounded, non-negative and non-increasing function with $V(1) = 0$ (a fully congested road results in zero velocity). In this paper, we only consider the normalized density meaning that $\rho$ is dimensionless.

Notice that in \eqref{equ:micro1}, the number of parameters increases with the number of vehicles, exposing the problem to the curse of dimensionality. In the FL1 model, the dynamics of the system are defined by how the position of vehicles changes over time. Through this assumption, the model is easy to implement in simulations, however, due to how humans perceive traffic in real life, it is more natural to consider how the velocity of vehicles evolves. This would imply to derive a second-order microscopic model that is not described in this paper but used in the simulation software.


%
For high-fidelity simulations, we use SUMO: an open-source software for microscopic traffic simulations. It is the short name for 'Simulation of Urban MObility' using second-order follow-the-leader model. It is used to construct a multimodal road network with various features such as junctions, route assignments, and car communication \citep{behrisch2011}. In this paper, SUMO is used to build a single-lane road with an attached driveway in the middle and a stochastic influx of vehicles. Moreover, a traffic light is placed on the highway, giving rise to congestion when vehicles from the driveway enter the highway. Nevertheless, for larger traffic systems, the microscopic model is, by design, limited computationally. As an alternative, we can use a macroscopic model that does not have the same restrictions.

\subsection{Macroscopic Model}
A macroscopic model uses three quantities to describe the temporal development of a traffic system, namely relative density $\rho$, average velocity $v$, and flow rate $q$, which all depend on time $t$ and space $x$. The relative density $\rho$ is the number of vehicles within a unit of road space normalized over the interval $[0,1]$, where $\rho=0$ is interpreted as an empty road and $\rho=1$ is a bumper-to-bumper road. Average velocity $v$ is the mean speed of vehicles on a road segment. Flow rate $q$ is the number of vehicles passing through a cross-section of the road during a time period. These three quantities are related through the hydrodynamic equation:
\begin{equation} \label{equ:hydro}
    q = \rho v.
\end{equation}
At the same time, the conservation of mass leads to the following equation:
\begin{equation} \label{equ:conservation}
    \rho_t + q_x = 0,
\end{equation}
Finally, we substitute \eqref{equ:hydro} for the flow rate $q$ in \eqref{equ:conservation} resulting in the continuity equation for traffic flow:
\begin{equation} \label{equ:continuity}
    \rho_t + (\rho v)_x = 0.
\end{equation}

The dynamics of Lighthill-Whitham Richards (LWR) model \citep{lighthill1955} are defined by
\begin{equation} \label{equ:macro1}
    \rho_t + (\rho V(\rho))_x = 0,
\end{equation}
where $V(\rho)$ is the equilibrium velocity which is analogous to the velocity function in the FL1 model. An important distinction to remark is that $v(t,x)$ is the average velocity at a position $x$ and a time $t$ while $V(\rho)$ is the velocity in an equilibrium state where the density is $\rho$. In the LWR model, it is assumed that $v(t,x) = V(\rho(t,x))$. Additionally, note that the velocity function in the FL1 model is not necessarily differentiable while the equilibrium velocity in LWR is required to be differentiable for a well-defined macroscopic model.

However, in this paper, we do not consider the system dynamics to be described by the continuity equation in \eqref{equ:macro1} due to the discontinuous nature of its solutions. Instead, as suggested by \cite{nelson2002}, we assume that the traffic flow is characterized by the following parabolic PDE:
\begin{equation} \label{equ:parabolic}
    \rho_t + (\rho V(\rho))_x = \gamma \rho_{xx},
\end{equation}
where $0 < \gamma \ll 1$ is a diffusion correction parameter. As $\gamma \rightarrow 0$, the solution of \eqref{equ:parabolic} converges, in a distributive sense, to the solution of the continuity equation in \eqref{equ:macro1}, which has been proven by \cite{hopf1950}. Equation \eqref{equ:parabolic} now exhibits a unique smooth solution which is favorable from a modeling and analysis point of view. Admittedly, equation \eqref{equ:parabolic} violates the anisotropic property of traffic flow but at the same time enables the possibility to model other traffic flow phenomena for which the LWR model has been proven to be deficient \citep{nelson2002}.

The main weakness of a first-order macroscopic model is the assumption that the system is always in an equilibrium state. For instance, the LWR model is unable to predict phenomena such as capacity drop and stop-and-go waves \citep{kontorinakia2017}. 
This paper is a preliminary study and therefore, we do not consider higher order traffic models. 

\subsection{Coupled Micro-Macro Model}
Both the microscopic and macroscopic models attempt to describe traffic flow, although, from separate standpoints they are not completely uncorrelated. \cite{colombo2014} showed that the density solution in the FL1 model $\rho_i$ in \eqref{equ:micro_density} converges to the density solution in the LWR model $\rho$, as the number of vehicles $n$ tends toward infinity while the total vehicle mass $m$ in the FL1 model remains constant. Specifically, the solution in the FL1 model describes the trajectory solution for the LWR model.
\begin{equation} \label{equ:micro_macro_relation}
    \lim_{n \rightarrow 0} \rho_i(t) = \rho(t, x_i)
\end{equation}
To exploit this relationship, we can construct a Coupled Micro-Macro (CMM) model, where we consider the parabolic PDE in \eqref{equ:parabolic}, as a substitute for the LWR model, and the equilibrium velocity given by the ODE in \eqref{equ:micro1} from the FL1 model.
\begin{equation} \label{equ:cmm}
    \begin{cases}
        \rho_t + (\rho V(\rho))_x = \gamma \rho_{xx} \\
        \dot{x}_i = V(\rho(x_i)),  & i=1,2,\dots,n.
    \end{cases}
\end{equation}
In this paper, we study the model in \eqref{equ:cmm}. The following subsection is devoted to the problem formulation in a more general context.

\subsection{Problem formulation}
The objective of this paper is to successfully reconstruct the true value of the desired quantity from a sparse data set.

\subsubsection{General problem}

Let $(\H, \| \cdot \|_{\H})$ be a semi-normed vector space, $z \in \H$ the desired quantity and $\H_c$ a subset in $\H$. We say that $\hat{z}$ is a \textbf{Partial-State Reconstruction} (PSR) of $z$ if
\begin{equation}
    \hat{z} \in \mathcal{R}, \quad 
    \mathcal{R} = \argmin_{\bar{z} \in \H_c} \| z - \bar{z} \|_{\H}^2,
\end{equation}
where $\mathcal{R}$ is the partial-reconstructed set.

This problem is a special case of the one introduced by \cite{barreauL4DC}. The authors point out there that if $z \in \H_c$ then $\hat{z} = z$ weakly. In contrast, if $z \notin \H_c$ and $H_c$ is a vector space, a PSR $\hat{z}$ will be the orthogonal projection of $z$ in $\H_c$, which is assumed to be non-empty.

\subsubsection{Traffic flow problem}

In this paper, we investigate the PSR of the density solution $\rho$ of the LWR model. However, the density $\rho$ is assumed to follow the parabolic equation resulting in the density reconstruction $\hat{\rho}$ being a solution of \eqref{equ:parabolic}. 
The work of \cite{hopf1950} proves the existence of a solution to \eqref{equ:parabolic} which enforces the set $\H_c$ to be non-empty.

The reconstruction is performed using density measurements in space and time. They are collected over a time window $[0,T]$ by Probe Vehicles (PV), with position $x_1, x_2, \dots, x_n$, that interact within a traffic stream obeying the ODE in \eqref{equ:cmm}. Specifically, we are interested in the PSR in the region between the first and last PV:
\begin{multline} \label{equ:reconstruction_domain}
    \O_r = \{ (t,x) \in \R^+ \times \R \ | \\ 
    x \in [x_1(t), x_n(t)], \  \forall t \in [0, T+\Delta T]\},
\end{multline}
which we consider to be the reconstruction domain. Note that we intend to predict the future density as well, where $\Delta T$ is the horizon of the future prediction.

In relation to the general problem, $\H = H^1(\O_r, [0,1])$ and $\| \cdot \|_\H$ is the traditional $\L_2$-norm. Additionally, $\H_c = C^\infty(\O_r, [0,1])$ and bounded due to the assumption in \eqref{equ:parabolic}. Note that the model only considers continuous functions while the density solution can exhibit discontinuous tendencies. Furthermore, measurements are only collected along PV trajectories, meaning that we are, effectively, attempting to solve the PDE in \eqref{equ:parabolic} without initial conditions. As a result, the reconstruction set $\mathcal{R}$ is most likely not a singleton. To solve this optimization problem, we utilize a neural network as described in the following section.

\section{Machine Learning} \label{section:machine_learning}

Machine learning methods are data-based approaches for pattern recognition and have been applied in numerous fields of research \citep{goodfellow2016}. Learning is a data-based approach to modeling non-linear, or perhaps even unknown, relationships and structures for real-world applications. The structure of a DNN was loosely inspired by the neural processes in the human brain. Although, DNN is an outdated model for the human brain it has proven to be useful when trying to predict complex patterns.


\subsection{Fundamentals of Neural Networks}
A DNN consists of processing units, referred to as neurons, which form fully connected sequential layers. Each neuron uses the output of all neurons from the previous layer as the input and returns a single value as the output. In general, layers can be divided into three categories by their functionality.

First, the input layer receives input values that are fed into the neural network where each neuron in the input layer corresponds to an input value. Second, the hidden layers process the input data, under the constraint of a network model. Last, the output layer uses the data from the hidden layers to return a final output. Analogous to the input layer, each node in the output layer corresponds to an output value.

Consider the $j$th neuron in layer $i$ with input vector $\textbf{x}_i \in \R^{n_{i-1}}$ where the elements $y_{i-1,k}$ in the input vector are the output values of the $k$th neuron in the previous layer, consisting of $n_{i-1}$ neurons. Moreover, the neuron has a corresponding weight vector $\textbf{w}_{i,j} \in \R^{n_{i-1}}$ and a bias $b_{i,j} \in \R$. The weight vector determines the influence of each element in the input vector and the bias acts as an offset value. The output $y_{i,j}$ of the neuron $(i,j)$ is a linear function of the input $\textbf{x}_i$ subject to an activation function $\mathcal{A}: \R \rightarrow \R$:
\begin{equation} \label{equ:neuron}
    y_{i,j} = \mathcal{A}(\textbf{w}_{i,j}^T \textbf{x}_i + b_{i,j}),
\end{equation}
where the activation function $\mathcal{A}$ is operating elementwise \citep{goodfellow2016}. We can now easily extend \eqref{equ:neuron} from individual neuron processes to the computation in entire layers:
\begin{equation} \label{equ:layer}
    \textbf{y}_i = \mathcal{A}\left(W_i \textbf{x}_i + \textbf{b}_i\right),
\end{equation}
where $W_i \in \R^{n_i \times n_{i-1}}$, in which each row is constructed from weights $\textbf{w}_{i, \cdot}$ corresponding to the neurons in layer $i$, and $\textbf{b}_i \in \R^{n_i}$ contains the biases for all neurons in layer $i$. The input for the following layer is then $\textbf{x}_{i+1} = \textbf{y}_i$.

\subsubsection{Network Size}

The network size of a neural network reflects the complexity of the model. Specifically, each hidden layer represents a layer of complexity with each neuron representing features within the layer. The state of different features is, ultimately, defined by a set of numbers that usually do not have an intuitive meaning. Therefore, neural networks are often viewed as a \textit{black box} tool.

One of the main trade-offs when constructing a neural network is deciding the size of the network. The network size determines the model which is used to describe the system. If the network size is insufficiently large, the neural network will be unable to fully model the system, meaning that the system complexity is greater than the model complexity. On the other hand, if the network size is overestimated, the neural network will model individual data points and features instead of the underlying structure of the system. Subsequently, the neural network will be prone to overfitting and generalize poorly over unknown data \citep{goodfellow2016}. To counteract this, we can introduce a regularization agent to restrict the region in parameter space. This is further discussed in Subsection~\ref{sebsection:pidl}.

\subsubsection{Estimation Cost Function}
A cost function is constructed to measure the inaccuracy of the model estimation in relation to the reference value. A common type of cost function for regression networks is a quadratic discriminant function, the most popular being the mean square error (MSE). Suppose we are given a training set $\{(t_i, x_i, \rho_i)\}$, the MSE estimation cost function is given by
\begin{equation} \label{equ:est_cost}
    \J_{\rm{est}} = MSE_{(\rho,\hat{\rho})} = \frac{1}{N_{\rm{est}}} \sum_{i=1}^{N_{\rm{est}}} | \rho_i - \hat{\rho}(t_i, x_i) |^2,
\end{equation}
where $\hat{z}$ is the model estimation of $z$. The cost is minimized by determining the optimal weights and biases. There exists a plethora of procedures which have been developed to perform parameter optimization.

The most used ones are based on the gradient descent algorithm. They are only capable of finding a global minimum if the cost function is strictly concave. Consequently, for an arbitrary smooth cost function, the gradient descent algorithm lacks mechanics to differentiate a local minimum from a global minimum. Nevertheless, this is not a crucial concern since it has been shown by \cite{choromanska2015} that the local minimums have values close to the global minimum if the complexity of the neural network is sufficiently large. However, the neural network might overfit the dataset. Thus, a regularization agent as described in the following subsection is needed to prevent over-fitting.


\subsection{Physics Informed Deep Learning} \label{sebsection:pidl}
Physics informed deep learning \citep{raissi2019, sirignano2018} is based on the same principle as the traditional deep learning but extended to respect the dynamics of the system. It is similar to a model-based approach.
\subsubsection{Physics Cost Function}
Specifically, we want to enforce the neural network to behave according to a physics property on chosen points in the system domain. Let $N_{\rm{phy}}$ physics points $\{(t_j, x_j)\}$, randomly selected through latin hypercube sampling on $\O_r$. Assume that the physical system obeys the relation $f(t,x,\rho) = 0$ on its domain. We can now define a physics cost function in \eqref{equ:phy_cost} using the physics property $f$.
\begin{equation} \label{equ:phy_cost}
    \J_{\rm{phy}} = MSE_{(0,f)} = \frac{1}{N_{\rm{phy}}} \sum_{j=1}^{N_{\rm{phy}}} | f(t_j, x_j, \hat{\rho}) |^2
\end{equation}
The total cost function $\J$ is constructed as the sum of the estimation cost and physics cost weighted by $\mu \in [0,1]$.
\begin{equation} \label{equ.tot_cost}
    \J = \mu \J_{\rm{est}} + (1 - \mu) \J_{\rm{phy}},
\end{equation}
The neural network is forced to make a trade-off between minimizing the estimation cost and the physics cost when optimizing the parameters $\Theta$ depending on the weight $\mu$. Qualitatively, using $\mu=0$ means that we are trying to solve a differential equation without initial or boundary conditions. It is unhelpful to consider this case since $z\equiv0$ is usually an optimal solution. In contrast, using $\mu=1$ is equivalent to a traditional neural network, which can result in a poor reconstruction in regions with insufficient measurements. Hence, an appropriate choice of $\mu$ is required for a desirable reconstruction.

From a machine learning perspective, the physics cost function acts as a regularization agent. The neural network aims at accurately estimating the measurement data under the constraint that the physics cost needs to remain low. Importantly, the physics cost is independent of the measurements and depends only on the parameter tensor $\Theta$. This additional cost function can be interpreted as introducing artificial measurements. Hence, it is an efficient solution to prevent overfitting. Moreover, in theory, it is a model-based approach which deters non-physical solutions and, therefore, it is more robust to noise in the measurements.


\subsection{Traffic flow problem}
For the traffic flow problem, we have $N_{\rm{est}}$ estimation points $\{(t_i,x_i,\rho_i)\}$, from PV following the ODE in \eqref{equ:cmm}. Consequently, the measurements are not well-spread in the reconstructed domain $\O_r$, indicating that using physics informed learning could be advantageous compared to traditional learning \citep{huang2020}. 
We assume the vehicle dynamics are described by \eqref{equ:cmm} such that
\begin{equation}
    f(t,x,\hat{\rho}) = \hat{\rho}_t(t,x) + (\hat{\rho} V(\hat{\rho}))_x(t,x) - \gamma \hat{\rho}_{xx}(t,x)
\end{equation}
is the physics property of the neural network. The total cost function for the traffic flow problem is constructed analogous to the example problem, where
\begin{align}
    \J_{\rm{est}} &= \frac{1}{N_{\rm{est}}} \sum_{i=1}^{N_{\rm{est}}} \sum_{j=1}^n | \rho(t_i, x_j(t_i)) - \hat{\rho}(t_i, x_j(t_i)) |^2, \\
    \J_{\rm{phy}} &= \frac{1}{N_{\rm{phy}}} \sum_{j=1}^{N_{\rm{phy}}} | f(t_j, x_j, \hat{\rho}) |^2,
\end{align}
\begin{equation}
    \J = \mu \J_{\rm{est}} + (1 - \mu) \J_{\rm{phy}}.
\end{equation}
Next, we present the reconstruction results from a finite difference scheme and a SUMO simulation using physics informed learning.

\section{Results} \label{section:results}

In this paper, we have attempted to reconstruct the traffic density from a finite difference scheme and a SUMO simulation. In the finite difference scheme, the vehicle dynamics were characterized by the equation in \eqref{equ:parabolic}, where the equilibrium velocity was based on the famous Greenshields equation \citep{greenshields1935}. In the SUMO simulation, vehicles follow a variation of a second-order microscopic model, which is capable of simulating numerous traffic flow phenomena \citep{kraub1998}. We have also investigated the reconstruction results for both simulations when noise was added to the measurements.

\subsection{Finite Difference Scheme}
\begin{figure}
    \centering
    \includegraphics[width=0.9\linewidth]{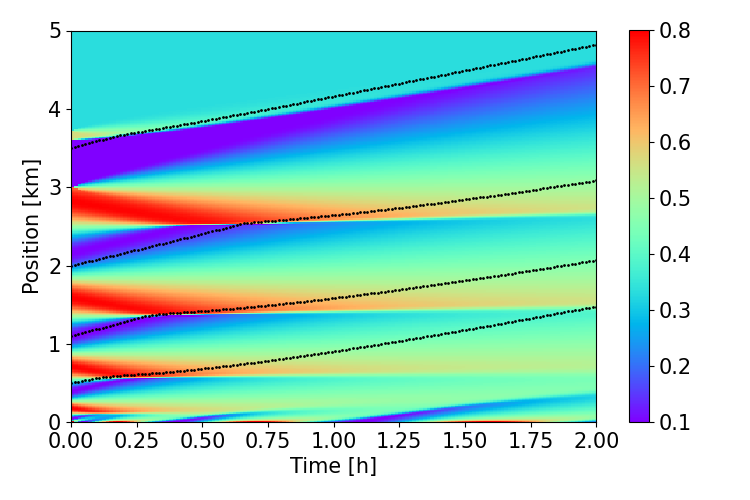}
    \caption{Spatiotemporal diagram of the true density for the finite difference scheme ($\gamma = 0.05$). The black points display the measurement points as well as the PV trajectories.}
    \label{fig:scheme_density_true}
\end{figure}

For the finite difference scheme, the vehicle velocity was determined by Greenshields equation:
\begin{equation}
    V(\rho) = v_f (1 - \rho),
\end{equation}
where $v_f$ is the vehicle velocity without interference, also called the free flow velocity. The true density solution is shown in Fig.~\ref{fig:scheme_density_true} using a spatiotemporal diagram. The colors display the density and the black points indicate the measurement points as well as the PV trajectories.

The reconstruction result when using only traditional learning is not shown here since a similar study has been conducted by \cite{huang2020}. The conclusion was that there is a perfect reconstruction of the density along the PV trajectories, however, it is far from the true density in between. Specifically, the neural network can predict a positive characteristic speed in congested regions. 

\begin{figure}
\centering
\subfloat[Noiseless measurements]{\includegraphics[width=0.9\linewidth]{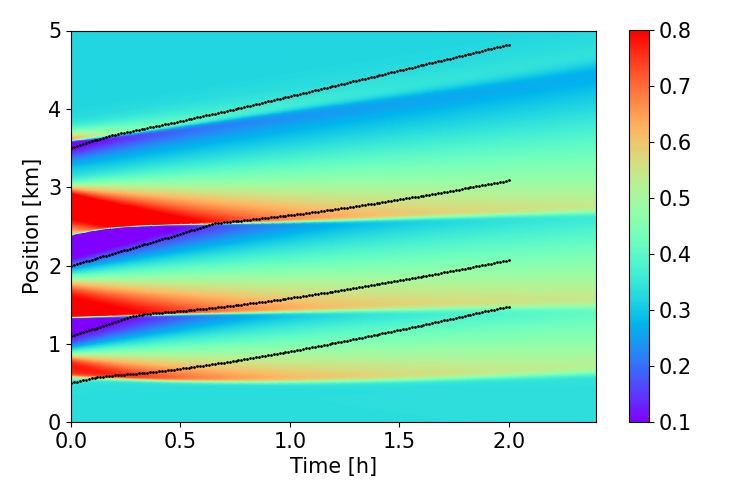}\label{subfig:scheme_lookahead}} \\
\vspace{-0.4cm}
\subfloat[With measurements noise]{\includegraphics[width=0.9\linewidth]{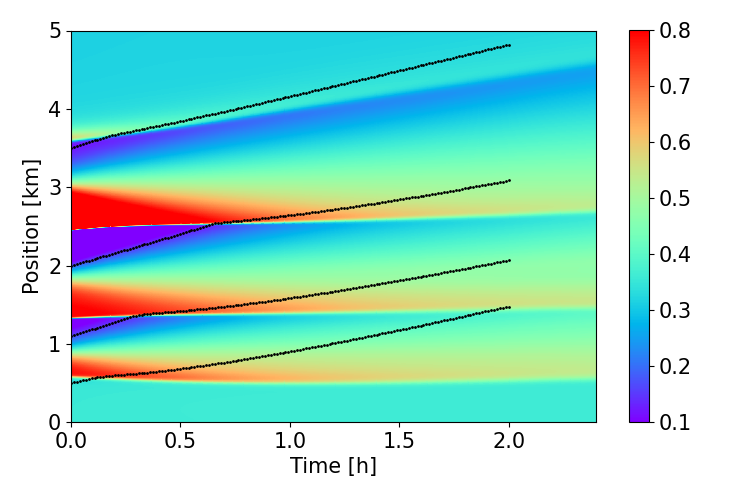} \label{subfig:scheme_lookahead_noise}}
\caption{Spatiotemporal diagrams of the traffic state reconstruction using physics informed deep learning. The black points are measurement points along the PV trajectories.}
\label{fig:scheme_reconstruction}
\end{figure}

For the physics informed reconstruction, in Fig.~\ref{subfig:scheme_lookahead}, we can observe that results are nearly indistinguishable from the true density for later times $t$, with an $\L_2$-error of $0.17$. Both cost functions are exhibiting low values. This showcases the ability of physics informed learning to compensate for insufficient measurement data. The neural network is also capable of performing accurate reconstruction on future density. The characteristics of the shock waves are preserved at all times in the reconstruction domain. Physics points are sampled throughout the entire reconstruction enforcing the physical property in the region $[T, T+\Delta T]$.

Additionally, we performed a reconstruction with the addition of an unbiased Gaussian noise in the measurements, shown in Fig.~\ref{subfig:scheme_lookahead_noise}. The standard deviation is $\sigma_\rho = 0.1$ (considerably higher than what would be expected in real-world scenarios). Surprisingly, the $\L_2$-error is $0.16$, so lower than the $\L_2$-error for the noiseless reconstruction. This is probably a coincidence, however, it shows the aptness of physics informed learning for unbiased noise rejection.

Admittedly, there are some discrepancies in both the noiseless reconstruction and the reconstruction with added noise, specifically at former times $t$. Nevertheless, this is to be expected since the neural network is unaware of the initial state of the system. From a control point of view, this is not an issue since the reconstruction of the past does not aid our ability to control the traffic flow.

\begin{figure}
    \centering
    \includegraphics[width=0.9\linewidth]{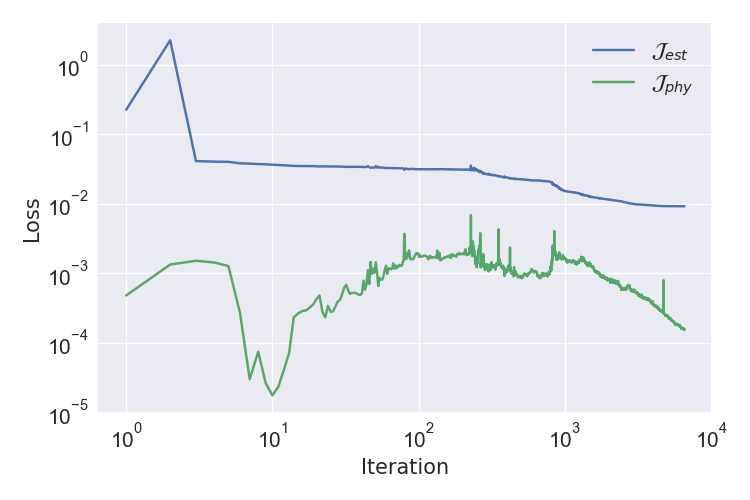}
    \caption{Values of the estimation cost (blue line) and the physics (green line) over all training iterations.}
    \label{subfig:loss}
\end{figure}

The training loss for the reconstruction with added noise is presented in Fig.~\ref{subfig:loss}. At the beginning of the training iterations, the estimation cost decreases while the physics cost increases, and for the later iterations, the network manages to reduce both costs simultaneously. It indicates the impact of the physics cost as a regularization agent. Specifically, the physics is never allowed to reach a value beyond the estimation cost.

\subsection{SUMO Simulation}

\begin{figure}
    \centering
    \includegraphics[width=0.9\linewidth]{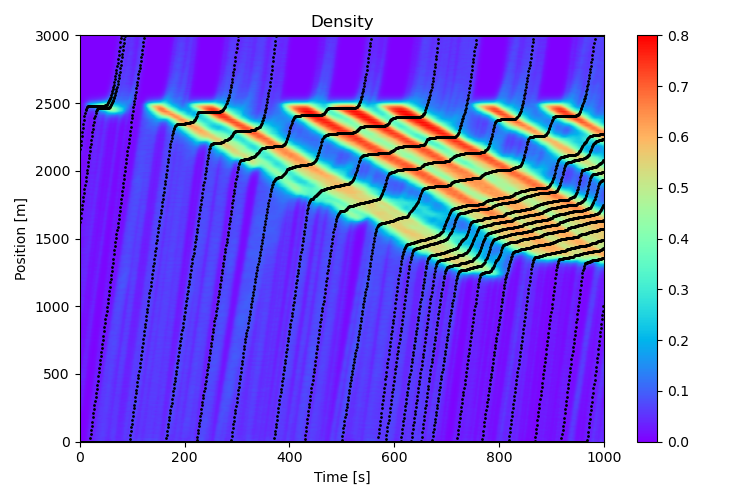}
    \caption{Spatiotemporal diagram of the true density for the SUMO simulation. The black points display the measurement points and the PV trajectories.}
    \label{fig:sumo_density_true}
\end{figure}

In the SUMO simulation, vehicles interact on a highway with a connected driveway and a traffic light as a mediator. Vehicles on the highway are stopped by the traffic light as vehicles from the driveway enter which induces congestion. The vehicles in the SUMO simulation follow a second-order microscopic model based on collision prevention. The true density is shown in Fig.~\ref{fig:sumo_density_true}.

\begin{figure}
\centering
\subfloat[Greenshields model]{\includegraphics[width=0.9\linewidth]{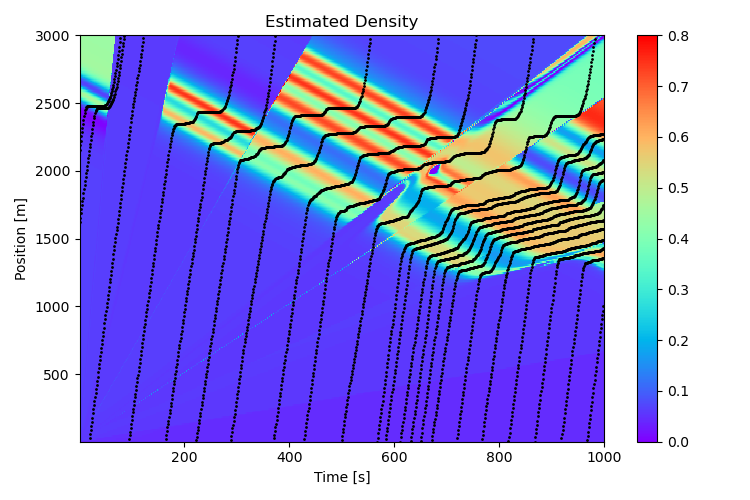} \label{subfig:sumo_greenshield}} \\
\subfloat[Model identification]{\includegraphics[width=0.9\linewidth]{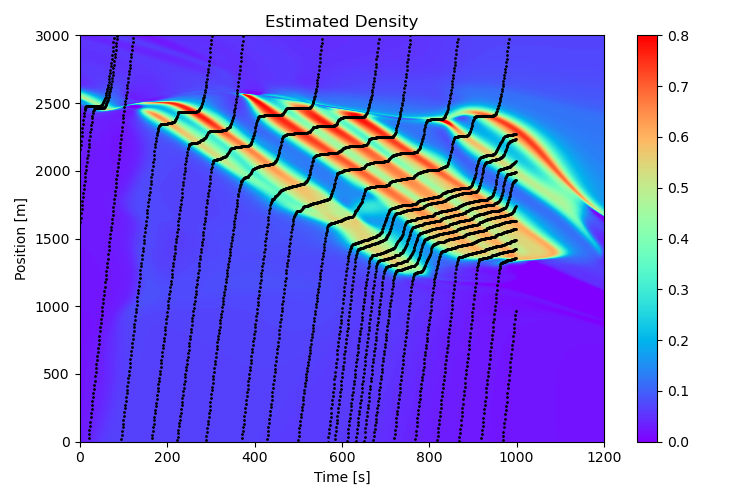} \label{subfig:sumo_pidl}}
\caption{Spatiotemporal diagrams of the physics informed reconstruction with different models. The black points are measurement points along the PV trajectories.}
\label{fig:sumo_reconstruction}
\end{figure}

The same method as in the previous subsection, where we assumed that vehicle dynamics are described by the Greenshields equation, result in the reconstruction in Fig.~\ref{subfig:sumo_greenshield}. We can observe that the neural network is incapable of reconstructing the true density. Greenshields equation assumes that there is a linear dependency between the density and the velocity whereas vehicles in the SUMO simulation exhibit non-linear behavior.

\begin{figure}
    \centering
    \includegraphics[width=0.9\linewidth]{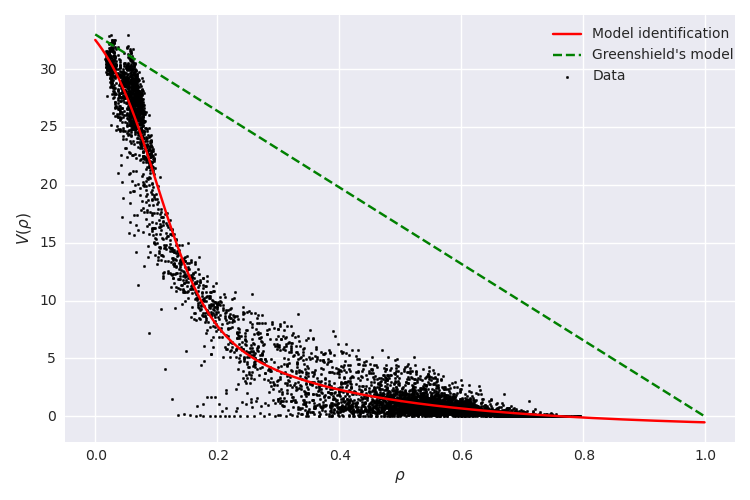}
    \caption{Equilibrium velocity resulting from Greenshields equation (green dotted line) and model identification (red line). The black points are empirical measurements.}
    \label{subfig:sumo_identification}
\end{figure}

Instead, we performed model identification using traditional learning to estimate the equilibrium velocity. Velocity measurements $\{v_i\}$ were constructed through numerical differentiation using the position and time measurements $\{t_i, x_i\}$. Fig.~\ref{subfig:sumo_identification} shows that the estimated velocity is significantly more accurate compared to the Greenshields equation. This difference is notable in the density reconstruction using model identification in Fig.~\ref{subfig:sumo_pidl}. Moreover, the neural network is also able to make future prediction, where the density characteristics are maintained.
The major discrepancy in the reconstruction is around $x=2500$. The reason for this is most likely due to the unlabelled presence of the traffic light.

The reconstruction with added Gaussian noise using $\sigma_\rho = 0.001$ is not shown here but it has a chaotic density which is considerably less accurate compared to the noiseless reconstruction. This may be due to insufficient model complexity but it could also be model instability. Further investigations are required to fully understand the reconstruction results with added noise.

\section{Conclusion} \label{section:conclusion}

In this paper, we have introduced a model-based approach for traffic state reconstruction using a neural network and presented the results for a finite difference scheme and a SUMO simulation. With the leverage of machine learning our method can distill the problem of identification, reconstruction, and prediction into one optimization problem. Additionally, the method is robust when introducing unbiased noise in density given a sufficient understanding of the system dynamics. Future work will be devoted to extend the method using 2nd order models and develop tools to counteract trajectory noise.

\bibliography{KEX_bib_example}

\end{document}